\documentclass[11pt]{article}
\usepackage{amssymb,amsthm,amsmath,tikz,pgf,mathtools}
\usepackage[lined,linesnumbered,ruled]{algorithm2e}
\usepackage{epsfig,amsmath, amsfonts,amstext, amsthm, latexsym, graphicx}
\usepackage{epstopdf}
\usepackage{graphicx}

\usepackage[top=1in, bottom=1in, left=1.6in, right=1.6in]{geometry}

 \newtheorem{fact}{Fact}

\begin{document}
\title{Some considerations on Collatz conjecture}

\author{Fabrizio Luccio\\
Department of Informatics, Universit\`{a} di Pisa\\
luccio@di.unipi.it
}


\maketitle
\begin{abstract}
Some simple facts are proved ruling the Collatz tree and the chains of vertices appearing in it, leading to the reduction of the number of significant elements appearing in the tree. Although the Collatz conjecture remains open, these fact may cast some light on the nature of the problem.
\end{abstract}


\section{Problem presentation}\label{sec:presentation}

The Collatz Conjecture (1937), also called $3x + 1$ Problem, concerns iteration of the map $T : Z^+ \rightarrow Z^+$ given by the two rules: 

\vspace{3mm}
{\bf R1} $\;T(x) =x/2$  if $x$ is even

\vspace{1mm}
{\bf R2} $\;T(x) =(3x + 1)/2$  if $x$ is odd

\vspace{3mm}
\noindent and asserts that each $n \geq 1$ has some iterate $T^{(h)}(n)=1$.
Hence each $n$ would eventually converge to the limit cycle $\;2 \rightarrow 1 \rightarrow 2 \dots$

\vspace{2mm}
There is no proof that the conjecture holds, that is for all the integers in $Z^+$ no cycle other than $\;2 \rightarrow 1 \rightarrow 2 \dots$ is formed, or in general that all integers iterate to 1, although experiments conducted on the integers up to $5 \times 2^{60}$ show that this is the case for them.
It is known, however, that some integers follow a long chain of transformations before getting to the limit cycle. For example $n=27$ takes 80 steps, climbing to 4616 before descending to 1, according to the sequence:

\vspace{3mm}
\noindent {\bf 27}, 41, 62, 31, 47, 71, 107, 161, 242, 121, 182, 91, 137, 206, 103, 155, 233, 350, 175, 263, 395, 1593, 890, 445, 668, 334, 167, 251, 377, 566, 283, 425, 638, 319, 479, 719, 1079, 1619, 2429, 3644, 1822, 911, 1367, 2051, 3077, {\bf 4616}, 2308, 1154, 577, 866, 433, 650, 325, 488, 244, 122, 61, 92, 46, 23, 35, 53, 80, 40, 20, 10, 5, 8, 4, 2, {\bf 1}

\begin{figure}[h]
\centering
\includegraphics[scale=0.4]{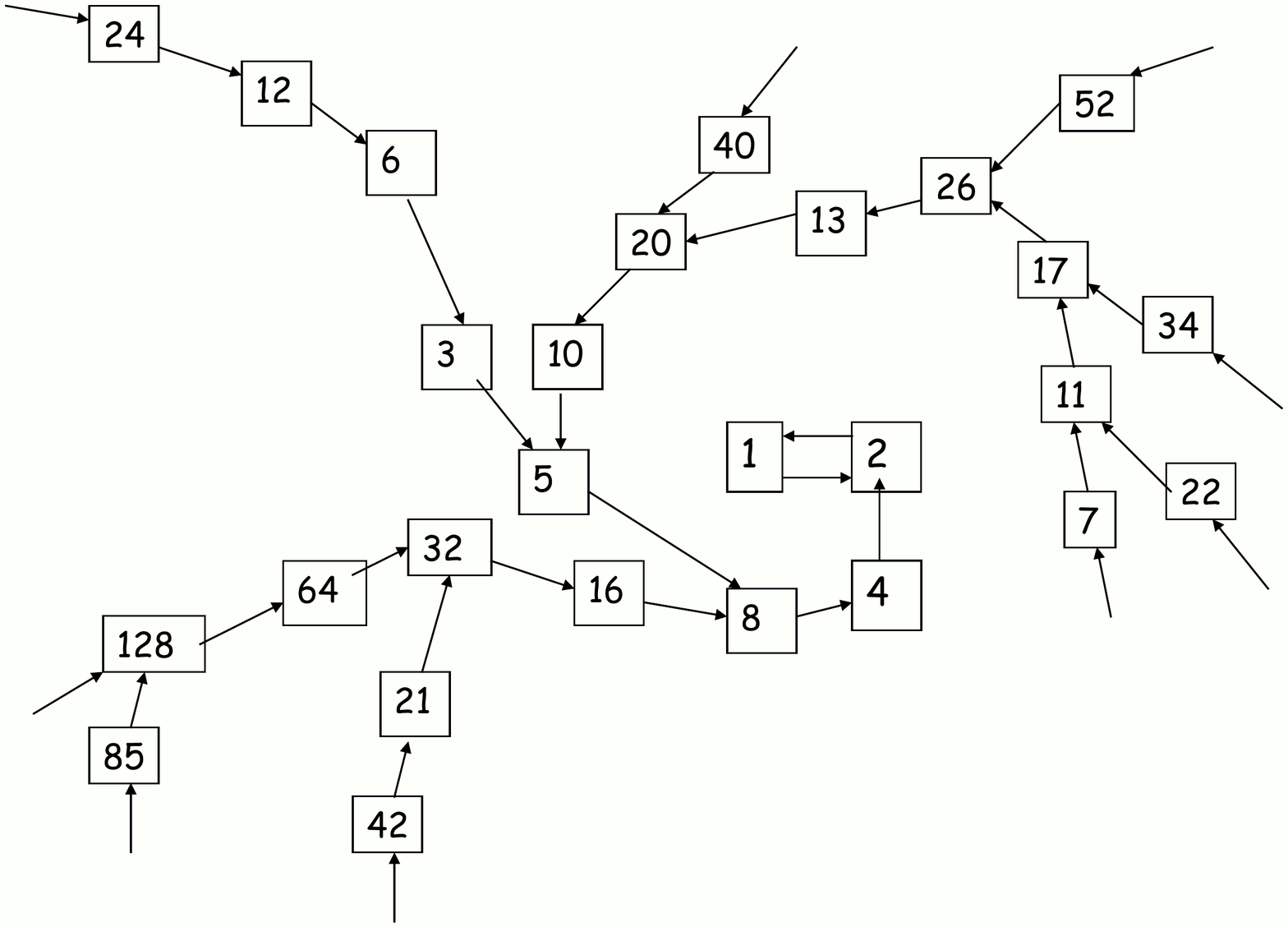}
\caption {\small Part of the Collatz tree.
\label{fig:tree}}
\end{figure}

\vspace{3mm}
The evolution of the integers induced by $T$ can be immediately represented as a directed graph that, according to the conjecture, is in fact a tree (called the {\em Collatz tree}) ending in the limit cycle at the root. A portion of this tree is shown in Figure~\ref{fig:tree}. 

\newpage
\section{Some properties of the Collatz tree}\label{sec:properties}

We now state some simple facts ruling the Collatz tree and the chains of vertices appearing in it. To the best of our knowledge these points have not been explicitly raised before. 

First note that each vertex $x$ of the graph may have one or two predecessors deriving from the application of {\bf R1} or {\bf R2} to a previous vertex. This leads to consider two inverses of $T$ indicated by $P_e(x)$ and $P_o(x)$, where the subscripts $e$ and $o$ indicate if the predecessor is even or odd, hence rule {\bf R1} or {\bf R2} has been applied. Clearly each vertex $x$ has an even predecessor $P_e(x)=2x$, and may have an odd predecessor $P_o(x)=(2x-1)/3$ (inverse of {\bf R2}) if this value is an odd integer. For example, for $x=5$ we have $P_e(5) = 10$ and $P_o(5)=3$, see Figure \ref{fig:tree}. 

In particular, for the non-negative integers consider the residue classes [0], [1], [2] modulo 3, respectively containing the integers $3k$, $3k+1$, $3k+2$ with $k\geq 0$. We have:

\begin{fact}
\label{fact:edges}

\vspace{1mm}
{\em (i)} each vertex $x\in [0]$ has only one predecessor $P_e(x)\in [0]$; 

\vspace{1mm}
{\em (ii)} each vertex $x\in [1]$ has only one predecessor $P_e(x)\in [2]$;

\vspace{1mm}
{\em (ii)} each vertex $x\in [2]$ has both predecessors $P_e(x)$ and $P_o(x)$, with$P_e(x)\in [1]$, 

\hspace{6mm}and $P_o(x)\in [0]$ if $(x-2)/3\in [1]$, $P_o(x)\in [1]$ if $(x-2)/3\in [0]$, $P_o(x)\in [2]$ 

\hspace{6.5mm}if $(x-2)/3\in [2]$.

\vspace{2mm}
\noindent {\bf Proof.} 
{\em Each vertex $x$ has an even predecessor $P_e(x)=2x$, and also has an odd predecessor $P_o(x)=(2x-1)/3$ if {\bf R2} applies to that predecessor. In case (i) we have $x=3k$ and $P_o(x)=(6k-1)/3$ and in case (ii) we have $x=3k+1$ and $P_o(x)=(6k+1)/3$, however, these two values are not integer and $P_o(x)$ does not exist. In case (iii) we have $x=3k+2$ and $P_o(x)=(6k+3)/3 = 2k+1$, an odd integer and evolves to $3k+2$ according to {\bf R2}. The residue classes of $P_o(x)$ according to the nature of $(x-2)/3$ are straightforwardly derived. }
\hfill $\Box$

\end{fact}

Referring to Figure \ref{fig:tree}, case (i) of Fact~\ref{fact:edges} applies to the backward chains $\;$3, 6, 12,  24, . . .   $\;$ and $\;$21, 42, . . .  $\;$ with $x=3k$. Case (ii) applies to vertices $\;$1, 4, 7, 10, 13, 16 . . .  $\;$ with $x=3k+1$. Case (iii) applies to vertices $\;$2, 5, 8,  11, 17, 32, . . .  $\;$ with $x=3k+2$.

\begin{fact}
\label{fact:3k}
For each vertex $x\in [0]$ we have:

\vspace{1mm}
{\em (i)} $T(x)\in [0]$ for $x/3$ even;

\vspace{1mm}
{\em (ii)} $T(x) \in [2]$ for $x/3$ odd.

\vspace{2mm}
\noindent {\bf Proof.} 
{\em Point (i) is immediate. For point (ii) we have $x=3k$ and $T(x)=(9k+1)/2 = 3h+r$, with $0\leq r \leq 2$, hence $h=(9k-2r+1)/6$. Since $k$ is odd, $h$ is integer only for $r=2$. }
\hfill $\Box$

\end{fact}

Referring to Figure \ref{fig:tree}, point (i) of Fact~\ref{fact:3k} applies to vertices 6 and 42, and point (ii) applies to vertices 3 and 21.

\begin{fact}
\label{fact:3k+1}
For each vertex $x\in [1]$ we have $T(x) \in [2]$.

\vspace{2mm}
\noindent {\bf Proof.} 
{\em $x=3k+1$. For $k$ even we have: $T(x)=(9k+4)/2 = 3h+r$, with $0\leq r \leq 2$, hence $h=(9k-2r+4)/6$. Since $k$ is even, $h$ is integer only for $r=2$. For $k$ odd we have $x=3(k-1)+4$ even, then $T(x)=3(k-1)/2 +2= 3h+2$ with $h=(k-1)/2$ integer.}
\hfill $\Box$

\end{fact}

Referring to Figure \ref{fig:tree}, Fact~\ref{fact:3k+1} applies to vertices 7 and 13 for $k$ even, and to vertices 4 and 10 for $k$ odd. In both cases $T(x)$ can be even or odd.

\begin{fact}
\label{fact:3k+2}
For each vertex $x\in [2]$ we have:

\vspace{1mm}
{\em (i)} $T(x) \in [1]$, for $x$ even;

\vspace{1mm}
{\em (ii)} $T(x) \in [2]$, for $x$ odd.

\vspace{2mm}
\noindent {\bf Proof.} 
{\em $x=3k+2$. For $x$ even we have $k$ even and $T(x)=3k/2+1$. For $x$ odd we have $k$ odd and $T(x)=(9k+7)/2 = 3h+r$, with $0\leq r \leq 2$, hence $h=(9k-2r+7)/6$. Since $k$ is odd, $h$ is integer only for $r=2$. }
\hfill $\Box$

\end{fact}

Referring to Figure \ref{fig:tree}, point (i) of Fact~\ref{fact:3k+2} applies to vertices 8 and 20, and point (ii) applies to vertices  5 and 11.

\begin{fact}
\label{fact:ciclo2}
The Collatz graph contains no cycle of length one and exactly only cycle of length two consisting of vertices 1 and 2. 

\vspace{2mm}
\noindent {\bf Proof.} 
{\em Obviously $T(x)\neq x$ hence no cycle of length one exists. In a cycle of length two we have $T^{(2)}(x)=x$ and both rules {\bf R1}, {\bf R2} must apply to avoid that $x$ iterates in two steps to a value that is certainly smaller or certainly larger than $x$. For $x$ even we have $T^{(2)}(x)=(3x/2+1)/2 = (3x+2)/4$, and the equality $(3x+2)/4=x$ implies $x=2$. For $x$ odd we have $T(x)=(3x+1)/2$. In this case, if $T(x)$ is even we have $T^{(2)}(x)=((3x+1)/2)/2 = (3x+1)/4$, and the equality $(3x+1)/4=x$ implies $x=1$;
if $T(x)$ is odd we have $T^{(2)}(x)>x$ since rule {\bf R2} applies twice.}
\hfill $\Box$

\end{fact}

Note that the existence of one or more finite cycles implies the existence of an infinite family of longer cycles, each made of the iteration of a shorter one like cycles 1-2-1-2, 1-2-1-2-1-2, etc. In the following, with the term cycle we refer to {\em simple} cycles where these iterations have been eliminated. 
Having said that, it is not known whether cycles exist in the Colatz graph, other than the one composed of vertices 1 and 2. As far as the existence of longer cycles, we can only state some partial results.
 
\begin{fact}
\label{fact:ciclo-x}
If the rules {\bf R1} and {\bf R2} are respectively applied $r_1$ and $r_2$ times in a cycle $C$, any vertex $x$ of $C$ must fulfil the equation:  $x=A/(2^{{r_1}+{r_2}}-3^{r_2})$, where $A$ is a positive integer and $r_1> 0.58\,r_2$.

\vspace{2mm}
\noindent {\bf Proof.} 
{\em Staring from vertex $x$, a sequence of $k=r_1+r_2$ applications of {\bf R1} and {\bf R2} must occur in $C$ until $x$ is reached again, that is $T^{(k)}(x)=x$. It is straightforward to verify that this leads to the equation: $(3^{r_2} x+ A)/ 2^{{r_1}+{r_2}} = x$, that is  $x=A/(2^{{r_1}+{r_2}} - 3^{r^2})$. As $x$ must be a positive integer we  have $2^{{r_1}+{r_2}} > 3^{r_2}$ hence $r_1> r_2 (log_2 3-1)$.
}
\hfill $\Box$

\end{fact}

We can easily prove that no cycle of length three exists. In fact, by Fact \ref{fact:ciclo-x} such a cycle should require applying {\bf R2} once and {\bf R1} twice, in one of the three orderings {\bf R1-R1-R2}, {\bf R1-R2-R1}, {\bf R2-R1-R1}, and in all cases the relation $T^{(3)}(x)=x$ cannot be satisfied. By the same reasoning, to have a cycle of length four one should apply {\bf R1} two or three times, and only the sequences {\bf R1-R2-R1-R2} and {\bf R2-R1-R2-R1} satisfy $T^{(4)}(x)=x$ with $x=2$ and $x=1$ respectively, corresponding to the non-simple cycles 2-1-2-1 and 1-2-1-2.
 
\begin{fact}
\label{fact:ciclo[0]}
No vertex $v\in [0]$ can belong to a finite cycle.

\vspace{2mm}
\noindent {\bf Proof.} 
{\em All the vertices in $[0]$ belong to some infinite chain $C$ of vertices $2^i\times 3k$ with $k\geq 1$ odd and $i\geq 0$, whose last vertex is $3k$ and each vertex has exactly one predecessor by Fact~\ref{fact:edges}(i). A finite cycle including one or more vertices of $C$ should include vertex $3k$ and procede through a chain of vertices not in $C$, one of which, say $x$,  has a successor $y$ in $C$ to form the cycle. Then $y$ should have two predecessors, one in $C$ and the other not in $C$, against Fact~\ref{fact:edges}(i)}.
\hfill $\Box$

\end{fact}

For example see the backwards chains 3, 6, 12, 24, . . .  and 21, 42, . . . in Figure~\ref{fig:tree}. 

\begin{fact}
\label{fact:ciclo[1]}
If a vertex $x\in [1]$ belongs to a finite cycle of length at least two, the existence of a cycle is preserved by eliminating $x$ from the graph and connecting $P_e(x)$ to $T(x)$ directly.

\vspace{2mm}
\noindent {\bf Proof.} 
{\em By Facts \ref{fact:edges} and \ref{fact:3k+1} each vertex $x\in [1]$ has exactly one predecessor and one successor in $[2]$.
 If $x$ belongs to a cycle $C$ of length $\ell\geq 2$, a new cycle of length $\ell-1$ arises connecting $P_e(x)$ to $T(x)$ directly}.
\hfill $\Box$

\end{fact}

For example, in Figure \ref{fig:tree} vertices 4 and 13 could be eliminated connecting 8 to 2, and 26 with 20, respectively.

As a consequence of Facts \ref{fact:ciclo[0]} and \ref{fact:ciclo[1]}, the existence of cycles can be studied maintaining only the vertices in [2] with the insertion of new edges as indicated in Fact \ref{fact:ciclo[1]}, and reformulating  
the Collatz conjecture as the iteration of the map $T' : [2] \rightarrow [2]$ given by the rules: 

\vspace{2mm}
{\bf R'1} $\;\;T'(x) =x/4$  if $x$ and $x/2$ are even

\vspace{1mm}
{\bf R'2} $\;\;T'(x) =(3x + 2)/4$  if $x$ is even and $x/2$ is odd

\vspace{1mm}
{\bf R'3} $\;\;T'(x) =(3x + 1)/2$  if $x$ is odd 


\vspace{2mm}
\noindent and asserting that each $n \in [2]$ has some iterate $T^{(h)}(n)=2$.

\begin{figure}[h]
\centering
\includegraphics[scale=0.35]{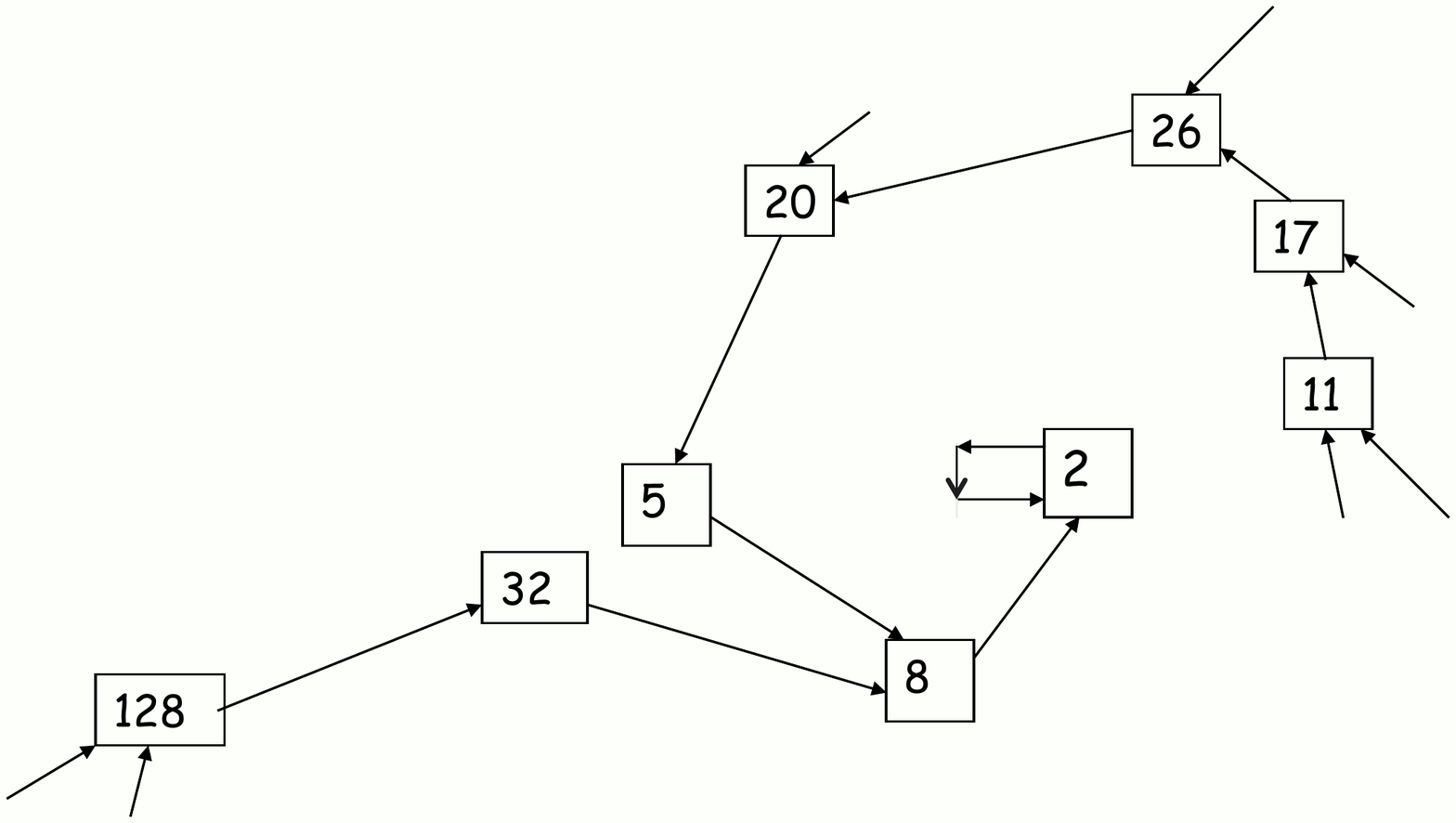}
\caption {\small The reduced Collatz tree.
\label{fig:reduced-tree}}
\end{figure}

\vspace{2mm}
It can be easily seen that if $x\in [2]$ also $T'(x)\in [2]$. In particular, compared with the original map $T$, rules {\bf R'1} and {\bf R'2} perform two steps while {\bf R'3} coincides with {\bf R2} of $T$, since in the original graph vertex $x$ is followed by another vertex in [2] (see Fact \ref{fact:3k+2}).
Figure \ref{fig:tree} then reduces to Figure \ref{fig:reduced-tree}. Note that {\bf R'1} is applied to vertices 8, 20, 32, and 128; {\bf R'2} is applied to vertices 2 and 26; {\bf R'3} is applied to vertices 5, 11, and 17.
The limit cycle 2 - 1 reduces to a cycle of length one containing the sole vertex 2.

\section{Bibliography}\label{sec:bibl}

A huge bibliography on Collatz conjecture exists. We just mention the excellent annotated list by J.C.Lagarias: ``The 3n+1 Problem: An Annotated Bibliography, II (2000-2009)'', arXiv:math/0608208v6 [math.NT] 12 Feb 2012.

\end{document}